\documentclass[12pt]{article}
\usepackage{amsmath}
\usepackage{graphicx}
\usepackage{CJK}
\usepackage{multirow}
\usepackage{makecell}

\usepackage{subfigure}
\usepackage[justification=centering]{caption}
\input{amssym.tex}

\def\bc{\begin{center}}       \def\ec{\end{center}}
\def\ba{\begin{array}}        \def\ea{\end{array}}
\def\be{\begin{equation}}     \def\ee{\end{equation}}
\def\bea{\begin{eqnarray}}    \def\eea{\end{eqnarray}}
\def\beaa{\begin{eqnarray*}}  \def\eeaa{\end{eqnarray*}}

\def\mathbb{\Bbb}

\begin{document}

\centerline {\bf \large On the number of limit cycles for Bogdanov-Takens system}
\vskip 0.1 true cm
\centerline {\bf \large under perturbations of piecewise smooth polynomials}

\vskip 0.3 true cm

\centerline{\bf  Jiaxin Wang, Liqin Zhao$^{*}$}
 \centerline{ School of Mathematical Sciences, Beijing Normal University,} \centerline{Laboratory of Mathematics and Complex Systems, Ministry of
Education,} \centerline{Beijing 100875, The People's Republic of China}

\footnotetext[1]{
This work was supported by NSFC(11671040)   \\ * Corresponding author.
E-mail:  zhaoliqin@bnu.edu.cn (L. Zhao).}

\vskip 0.5 true cm
\noindent{\bf  Abstract} In this paper, we study the bifurcate of limit cycles for Bogdanov-Takens system($\dot{x}=y$, $\dot{y}=-x+x^{2}$) under perturbations of piecewise smooth polynomials of degree $2$ and $n$ respectively. We bound the number of zeros of first order Melnikov function which controls the number of limit cycles bifurcating from the center. It is proved that the upper bounds of the number of limit cycles with switching curve $x=y^{2m}$($m$ is a positive integral) are $(39m+36)n+77m+21(m\geq 2)$ and $50n+52(m=1)$ (taking into account the multiplicity). The upper bounds number of limit cycles with switching lines $x=0$ and $y=0$ are 11 (taking into account the multiplicity) and it can be reached.

\vskip 0.1 true cm
\noindent{\bf Key Words} {limit cycle; Abelian integral; bifurcation. }

\vskip 0.5 true cm
\centerline{\bf{ $\S$1}. Introduction and the main results}
\vskip 0.3 true cm

\noindent The determination of limit cycles is one important problem in the qualitative theory of  planar differential systems. Stimulated by non-smooth phenomena in the real world such as control system [1], economics [12], nonlinear oscillations [20], and biology [5], [14], the investigation of limit cycles for piecewise smooth differential systems has attracted many attentions.

Many scholars have studied the number of limit cycles for piecewise smooth differential systems. In [11] and [17], the authors studied the expressions of the first order Melnikov function for the  piecewise Hamiltonian systems under the piecewise perturbations.

For the piecewise smooth differential systems with the switching lines, there have been a lot of results for example [4,16,18,23,26]. Nowadays many scholars begin to pay attention to the study of piecewise smooth differential systems with the nonlinear switching curves (see [2,6,19,21,22,24,27-29]).

It is shown by Horozov and Iliev in [9] that any cubic Hamiltonian can be transformed into the following normal form
$$H(x,y)=\frac{1}{2}(x^{2}+y^{2})-\frac{1}{3}x^{3}+axy^{2}+\frac{1}{3}by^{3},$$
where $a$, $b$ are parameters lying in the region
$$G=\left\{(a,b):-\frac{1}{2}\leq a\leq1,0\leq b\leq (1-a)(1+2a)^{1/2}\right\}.$$
Moreover, their respective vector fields $X_H$ are degenerate if $(a,b)\in \partial G$. If $(a,b)\in \partial G$, then in suitable coordinates (see [13]) all respective basic dynamics of $X_H$ can be classified into eight types that contain Bogdanov-Takens system with the fisrt integral
$$H(x,y)=\frac{1}{2}x^2+\frac{1}{2}y^{2}-\frac{1}{3}x^{3}=h,~h\in(0,\frac{1}{6}).\eqno(1.1)$$

Motivated by [23-25], in the present paper, we study the upper bounds of the number of limit cycles bifurcating from the period annuluses of Bogdanov-Takens system when it are perturbed inside any discontinuous polynomial differential systems. Concretely, we consider the following systems ($0<|\varepsilon|\ll 1$)
 $$
\left(
  \begin{array}{c}
          \dot{x}\\
          \dot{y}
   \end{array}
   \right)=\begin{cases}
   \left(
    \begin{array}{c}
        \quad y+\varepsilon p^+(x,y)\\
        -x+x^{2}+\varepsilon q^+(x,y)
   \end{array}
   \right), \quad x>y^{2m},\\
   \,
   \left(
    \begin{array}{c}
       \quad y+\varepsilon p^-(x,y)\\
       -x+x^{2}+\varepsilon q^-(x,y)
   \end{array}
   \right), \quad x<y^{2m},\\
   \end{cases}
   \eqno(1.2)$$
where $m$ is a positive integral and
$$p^{\pm}(x,y)=\sum\limits_{i+j=0}^na_{i,j}^{\pm}x^{i}y^{j},
~~q^{\pm}(x,y)=\sum\limits_{i+j=0}^nb_{i,j}^{\pm}x^{i}y^{j},$$
and
$$
\left(
  \begin{array}{c}
          \dot{x}\\
          \dot{y}
   \end{array}
   \right)=\begin{cases}
   \left(
    \begin{array}{c}
        \quad y+\varepsilon p^+(x,y)\\
        -x+x^{2}+\varepsilon q^+(x,y)
   \end{array}
   \right), \quad x>0, y>0,\\
   \,
   \left(
    \begin{array}{c}
       \quad y+\varepsilon \tilde{p}^+(x,y)\\
       -x+x^{2}+\varepsilon \tilde{q}^+(x,y)
   \end{array}
   \right), \quad x>0, y<0,\\
   \
   \left(
    \begin{array}{c}
       \quad y+\varepsilon \tilde{p}^-(x,y)\\
       -x+x^{2}+\varepsilon \tilde{q}^-(x,y)
   \end{array}
   \right), \quad x<0, y<0,\\
   \
   \left(
    \begin{array}{c}
       \quad y+\varepsilon p^-(x,y)\\
       -x+x^{2}+\varepsilon q^-(x,y)
   \end{array}
   \right), \quad x<0, y>0,\\
   \end{cases}
   \eqno(1.3)$$
where
$$p^{\pm}(x,y)=\sum\limits_{i+j=0}^2p_{i,j}^{\pm}x^{i}y^{j},
~~q^{\pm}(x,y)=\sum\limits_{i+j=0}^2q_{i,j}^{\pm}x^{i}y^{j},$$
$$\tilde{p}^{\pm}(x,y)=\sum\limits_{i+j=0}^2\tilde{p}_{i,j}^{\pm}x^{i}y^{j},
~~\tilde{q}^{\pm}(x,y)=\sum\limits_{i+j=0}^2\tilde{q}_{i,j}^{\pm}x^{i}y^{j}.$$
Let $H(n)$ be the number of limit cycles for system (1.2) and (1.3) bifurcating from the period annulus (taking into account the multiplicity). The main results are as follows.

\vskip 0.2 true cm

\noindent{\bf Theorem 1.1.} Consider system (1.2), by using the first order of Melnikov function in $\varepsilon$, the upper bounds of the number of limit cycles (taking into account the multiplicity) bifurcating from period annuli are

\noindent(i) If the switching curve is $y=x^{2}$, then $H(n)\leq 50n+52$.

\noindent(ii) If the switching curve is $y=x^{2m}(m\geq 2)$, then $H(n)\leq (39m+36)n+77m+21$.

\vskip 0.2 true cm

\noindent{\bf Theorem 1.2.} Consider system (1.3), by using the first order of Melnikov function in $\varepsilon$, the upper bounds of the number of limit cycles (taking into account the multiplicity) bifurcating from period annuli are 11, and the upper bounds can be reached for some $p_{i,j}^{\pm}(\tilde{p}_{i,j}^{\pm})$ and $q_{i,j}^{\pm}(\tilde{q}_{i,j}^{\pm})~(i,j=0,1,2)$.

\vskip 0.2 true cm

\noindent{\bf Remark 1.3.} B. Li et al. [15] considered respectively systems $(1.2)_{\varepsilon=0}$ under continuous perturbations of arbitrary polynomials with degree $n$. It is proved that for perturbed system $(1.2)_{\varepsilon=0}$, the exactly upper bound of the first order Melnikov function (Abelian integral) is $n-1$, and the exactly upper bound of the second order Melnikov function is $2n-2$ ($n$ is even) or $2n-3$ ($n$ is odd) when the first order Melnikov function vanishes. S. Sui et al. [23] and W. Cui et al. [3] considered respectively systems (1.1) under discontinuous perturbations of arbitrary polynomials with degree $n$ with switching line $y=0$ and $x=0$. It is proved that the upper bound of number of the isolated zeros of Abelian integrals for perturbed Bogdanov-Takens system are $12n+\left[\frac{n}{2}\right]+5$ and $16n+\left[\frac{n}{2}\right]-10$ respectively.

\vskip 0.5 true cm
\centerline{\bf{ $\S$2}. Preliminaries}
\vskip 0.3 true cm

Next, we shall introduce the first order Melnikov function of discontinuous differential systems. For $0<|\epsilon|\ll 1$, we consider the following Near-Hamilton system:
$$(\dot{x},\ \dot{y})=\begin{cases}
       (H^{+}_{y}(x,y)+\epsilon p^+(x,y),-H^{+}_{x}(x,y)+\epsilon q^+(x,y)),\ \ x\geq \psi(y),\\
       (H^{-}_{y}(x,y)+\epsilon p^-(x,y),-H^{-}_{x}(x,y)+\epsilon q^-(x,y)),\ \ x<\psi(y),\end{cases}\eqno(2.1)_\epsilon$$
where $\psi(x)$ is analytic with $\psi(0)=0$, and $p^\pm(x,y)$ and $q^\pm(x,y)$ are polynomials with degree $n$. System $(2.1)_\epsilon$ has two sub-systems:
$$\left\{{\begin{aligned}
\dot{x}&=~~~H^{+}_{y}(x,y)+\epsilon p^{+}(x,y),\\
\dot{y}&=-H^{+}_{x}(x,y)+\epsilon q^{+}(x,y),
\end{aligned}}~~~~~~x\geq \psi(y),\right.\eqno(2.2)$$
and
$$\left\{{\begin{aligned}
\dot{x}&=~~~H^{-}_{y}(x,y)+\epsilon p^{-}(x,y),\\
\dot{y}&=-H^{-}_{x}(x,y)+\epsilon q^{-}(x,y),
\end{aligned}}~~~~~~x<\psi(y).\right.\eqno(2.3)$$
Suppose that $(2.1)_{\epsilon=0}$ has a family of periodic orbits around the origin and satisfies the following assumptions.

\vskip 0.2 true cm

{\bf Assumption (I).} There exists an open interval $\Sigma$ such that for each $h\in\Sigma$, there are two points $A(h)$ and $B(h)$ on the curve $x=\psi(y)$ with $A(h)=(\psi(a(h)),a(h))$ and $B(h)=(\psi(b(h)),b(h))$ satisfying
$$H^{+}(A(h))=H^{+}(B(h))=h,~~H^{-}(A(h))=H^{-}(B(h)),~~a(h)<0<b(h).$$

{\bf Assumption (II).} The subsystem $(2.2)_{\epsilon=0}$ has an orbital arc $L_{h}^{+}$ starting from $A(h)$ and ending at $B(h)$ defined by $H^{+}(x,y)=h$ ($x\geq \psi(y)$). The subsystem $(2.3)_{\epsilon=0}$ has an orbital arc $L_{h}^{-}$ starting from $B(h)$ and ending at $A(h)$ defined by $H^{-}(x,y)=\tilde{h}$($x<\psi(y)$).

{\bf Assumption (III).} For each $h\in\Sigma$,
$$H_{x}^{\pm}(x,y)\psi^{'}(y)+H_{y}^{\pm}(x,y)\neq0~~~\text{at points $A(h)$ and $B(h)$}.$$
This means that the orbital arcs $L_{h}^{\pm}$ are not tangent to curve $x=\psi(y)$ at points $A(h)$ and $B(h)$.

\vskip 0.2 true cm

Under the assumptions {\bf Assumption (I),(II)} and {\bf (III)}, system $(2.1)_{\epsilon=0}$ has a family of non-smooth periodic orbits $L_{h}=L_{h}^{+}\cup L_{h}^{-}(h\in \Sigma)$. For definiteness, we assume that the orbits $L_{h}$ orientate clockwise. For $x=\psi(x)$, the authors [24] established the bifurcation function $F(h,\epsilon)$ for $(2.1)_\epsilon$. Let $F(h,0)=M(h)$. In [24], the authors obtained the following results.

\vskip 0.2 true cm

\noindent{\bf Lemma 2.1.}[24] Under the assumptions {\bf Assumption (I),(II)} and {\bf (III)}, we have

{\bf (i)}~If $M(h)$ has $k$ zeros in $h$ on the interval $\Sigma$ with each having an odd multiplicity, then $(2.1)_\epsilon$ has at least $k$ limit cycles bifurcating from the period annulus for $0<\left|\epsilon\right|\ll1$.

{\bf (ii)}~If $M(h)$ has at most $k$ zeros in $h$ on the interval $\Sigma$, taking into account the multiplicity, then there exist at most $k$ limit cycles of $(2.1)_\epsilon$ bifurcating from the period annulus.

{\bf (iii)} The first order Melnikov function $M(h)$ of system $(2.1)_\epsilon$ can be expressed as
$$M(h)=\int_{L_{h}^{+}}q^{+}dx-p^{+}dy+\frac{H_{x}^{+}(A)\psi^{'}(a(h))
+H_{y}^{+}(A)}{H_{x}^{-}(A)\psi^{'}(a(h))
+H_{y}^{-}(A)}\int_{L_{h}^{-}}q^{-}dx-p^{-}dy.\eqno(2.4)$$

\noindent{\bf Definition 2.2.}[7] We say that $\mathbb{V}$ is a {\bf Chebyshev space}, provided that each non-zero function in $\mathbb{V}$ has at most $dim(\mathbb{V})-1$ zeros, counted with multiplicity.

Let $\mathbb{S}$ be the solution space of a second order linear analytic differential equation
$$
x^{''}+a_{1}(t)x^{'}+a_{2}(t)x=0
\eqno(2.5)$$
on an open interval $\mathbb{I}$.

\noindent{\bf Lemma 2.3.}[7] The solution space $\mathbb{S}$ of $(2.5)$ is a Chebyshev space of the interval $\mathbb{I}$ if and only if there exists a nowhere vanishing solution $x_{0}(t)\in \mathbb{S}$($x_{0}(t)\neq0, \forall t\in \mathbb{I}$).

\noindent{\bf Lemma 2.4.}[7] Suppose the solution space of the homogeneous equation $(2.5)$ is a Chebyshev space and let $R(t)$ be an analytic function on $\mathbb{I}$ having $l$ zeros (counted with multiplicity). Then every solution $x(t)$ of the non-homogeneous equation
$$x^{''}+a_{1}(t)x^{'}+a_{2}(t)x=R(t)$$
has at most $l+2$ zeros on $\mathbb{I}$.

In this section we first introduce some results for determining the numbers of isolated zeros of a function.

\vskip 0.2 true cm

\noindent{\bf Definition 2.5.}[8] Let ${\cal F}=(f_0(x), f_1(x),...,f_{n}(x))$ be an ordered set of $C^\infty$ functions on an open interval $J\subset\mathbb{R}$. The ordered set $\left(f_0(x),f_1(x),...,f_{n}(x)\right)$ is said to be an ECT-system on $J$ if, for all $k=1,2,...,n,n+1$, any nontrivial linear combination
$$\alpha_0f_0(x)+\alpha_1f_1(x)+\cdots+\alpha_{k-1}f_{k-1}(x)$$
has at most $k-1$ isolated zeros on $J$ counted with multiplicities.

\vskip 0.2 true cm

\noindent{\bf Lemma 2.6.}[8] The ordered set $(f_0(x),f_1(x),...,f_{n}(x))$ is an ECT-system on $J$ if and only if, for each $k=1,2,...,n+1,$
$$W(f_0,f_1,...,f_{k-1})\neq0$$
for all $x\in J,$
where $W(f_0,f_1,...,f_{k-1})$ is the Wronskian of functions $f_0(x),f_1(x),...,f_{k-1}(x).$

\vskip 0.5 true cm
\centerline{\bf{ $\S$3}. Proof of Theorem 1.1}
\vskip 0.3 true cm

This section we consider system (1.2). System $(1.2)_{\varepsilon=0}$ has two singular points. There are a center $O(0,0)$ corresponding to $h=0$, a saddle $S(1,0)$ corresponding to $h=\frac{1}{6}$. For $h\in(0,\frac{1}{6})$, by Lemma 2.1, we have
$$M(h)=\int_{L_{h}^{+}}q^{+}(x,y)dx-p^{+}(x,y)dy+\int_{L_{h}^{-}}q^{-}(x,y)-p^{-}(x,y)dy, \eqno(3.1)$$
where
$$L_{h}^{+}(L_{h}^{-})=\{(x,y)|H(x,y)=h,x>y^{2}(x<y^{2m})\}.$$
Suppose that $H(x,y)=h$ and $y=0$ intersects at $(\bar{x}(h),0)$ and $(\tilde{x}(h),0)$ $(\bar{x}<0<\tilde{x})$. $x=y^{2m}$ and $H(x,y)=h$ intersects at points $A(u(h)^{2m},u(h))$ and $B(u(h)^{2m},-u(h))$. In the following we denote $u(h)$, $J_{i,j}(h)$ and $I_{i,j}(h)$ as $u$, $J_{i,j}$ and $I_{i,j}$. Hence we have
$$\frac{1}{2}u^{2}+\frac{1}{2}u^{4m}-\frac{1}{3}u^{6m}=h. \eqno(3.2)$$
For $i,j\geq 0$, let
$$J_{i,j}(h)=\int_{{L}_{h}^{+}}x^{i}y^{j}dx,~~I_{i,j}(h)=\int_{L_{h}^{-}}x^{i}y^{j}dx.$$
Without loss of generality, we only consider the case of $m=1$. For $m\geq 2$ it can be proved similarly. For polynomial $f(h)$, we denote it as $f(u)$ if we substitute (3.2) into $f(h)$.
\vskip 0.3 true cm

\noindent{\bf Lemma 3.1.} Consider system (1.2) for $h\in(0,\frac{1}{6})$. Then $M(h)$ can be expressed as
$$M(h)=\alpha(h)J_{0,1}+\beta(h)J_{1,1}+\gamma(h)I_{0,1}+\eta(h)I_{1,1}+\Phi(u), \eqno(3.3)$$
where ${\rm deg}\alpha(h),~\gamma(h)\leq \left[\frac{n-1}{2}\right]$, ${\rm deg}\beta(h),~\eta(h)\leq \left[\frac{n}{2}\right]-1$, and $\Phi(u)=\sum\limits_{i=0}^{3\left[\frac{n-1}{2}\right]+1}c_{i}u^{2i+1}$($c_{i}$ are constants).

\vskip 0.1 true cm
\noindent{\bf Proof.} First we assert that
$$M(h)=\sum_{i+2j+1=0}^{n}\rho^{+}_{i,2j+1}J_{i,2j+1}+\sum_{i+2j+1=0}^{n}\rho_{i,2j+1}^{-}I_{i,2j+1}+\sum_{i+2j=0}^{n}(a_{i,2j}^{+}-a_{i,2j}^{-})\frac{2}{2j+1}u^{2i+2j+1},\eqno(3.4)$$
Using Green's formula we have
$$\begin{aligned}\int_{L_{h}^{+}}x^{i}y^{j}dy=&-\frac{i}{j+1}\left(\int_{L_{h}^{+}}x^{i-1}y^{j+1}dx+\int_{\widehat{BOA}}x^{i-1}y^{j+1}dx\right)-\int_{\widehat{BOA}}x^{i}y^{j}dy\\
=&-\frac{i}{j+1}J_{i-1,j+1}-\frac{1+(-1)^{j}}{j+1}u^{2i+j+1}.
\end{aligned}$$
Similarly, we can obtain
$$\int_{L_{h}^{-}}x^{i}y^{j}dy=-\frac{i}{j+1}J_{i-1,j+1}+\frac{1+(-1)^{j}}{j+1}u^{2i+j+1}.$$
Thus,
$$
M(h)=\sum_{i+j=0}^{n}\rho^{+}_{i,j}J_{i,j}+\sum_{i+j=0}^{n}\rho_{i,j}^{-}I_{i,j}+\sum_{i+j=0}^{n}(a_{i,j}^{+}-a_{i,j}^{-})\frac{1+(-1)^{j}}{j+1}u^{2i+j+1},$$
where $u=u(h)$, $\rho_{i,j}^{\pm}=b_{i,j}^{\pm}+\frac{i+1}{j}a_{i+1,j-1}^{\pm}~(j\geq 1)$ and $\rho_{i,0}^{\pm}=b_{i,0}^{\pm}$.

Noticing the symmetry of $H(x,y)$, we have $J_{i,2j}=I_{i,2j}=0$. Therefore we can get (3.4). Next, differentiating (1.1) with respect to $x$, we obtain
$$x+y\frac{\partial{y}}{\partial{x}}-x^{2}=0.\eqno(3.5)$$
Multiplying (1.1) and (3.5) by $x^{i}y^{j}dx$ and integrating over $L_{h}^{+}$, we have
$$J_{i+1,j}-\frac{i}{j+1}J_{i-1,j+2}-J_{i+2,j}-\frac{1+(-1)^{j+1}}{j+2}u^{2i+j+2}=0,\eqno(3.6)$$
$$\frac{1}{2}J_{i+2,j}+\frac{1}{2}J_{i,j+2}-\frac{1}{3}J_{i+3,j}=hJ_{i,j}.\eqno(3.7)$$
Elementary manipulations reduce Eqs. (3.6) and (3.7) to
$$J_{i,j}=\frac{6j}{2i+3j+2}hJ_{i,j-2}-\frac{j}{2i+3j+2}J_{i+2,j-2}-\frac{2(1+(-1)^{j-1})}{2i+3j+2}u^{2i+j+2},\eqno(3.8)$$
$$J_{i,j}=-\frac{6(i-2)}{2i+3j+2}hJ_{i-3,j}+\frac{3(i+j)}{2i+3j+2}J_{i-1,j}-\frac{3(1+(-1)^{j+1})}{2i+3j+2}u^{2i+j-2}.\eqno(3.9)$$
Similarly, we have
$$I_{i,j}=\frac{6j}{2i+3j+2}hI_{i,j-2}-\frac{j}{2i+3j+2}I_{i+2,j-2}+\frac{2(1+(-1)^{j-1})}{2i+3j+2}u^{2i+j+2},\eqno(3.10)$$
$$I_{i,j}=-\frac{6(i-2)}{2i+3j+2}hI_{i-3,j}+\frac{3(i+j)}{2i+3j+2}I_{i-1,j}+\frac{3(1+(-1)^{j+1})}{2i+3j+2}u^{2i+j-2}.\eqno(3.11)$$
By (3.6), let $i=0$, $j=1$, we can obtain
$$J_{2,1}=J_{1,1}-\frac{2}{3}u^{3},\eqno(3.12)$$
Similarly, we can obtain
$$I_{2,1}=I_{1,1}+\frac{2}{3}u^{3}.\eqno(3.13)$$
Then using (3.4) and (3.8)--(3.13) we can obtain (3.3). This ends the proof.

\vskip 0.3 true cm

\noindent{\bf Lemma 3.2.} Let ${\bf V_{1}(h)}=(J_{0,1},J_{1,1})^{T}$, ${\bf V_{2}(h)}=(I_{0,1},I_{1,1})^{T}$, and $\sigma(u)=1/(1+2u^{2}-2u^{4})$. Then the vector functions ${\bf V_{1}(h)}$ and ${\bf V_{2}(h)}$ satisfy respectively the following Picard-fuchs equations:
$${\bf V_{1}(h)}=(B_{1}h+C_{1}){\bf V_{1}^{'}(h)}+W(u),\eqno(3.14)$$
$${\bf V_{2}(h)}=(B_{1}h+C_{1}){\bf V_{2}^{'}(h)}-W(u),\eqno(3.15)$$
where
$$B_{1}h+C_{1}=\left(\begin{matrix}
             &\frac{6}{5}h~~&-\frac{1}{5}\\
             &\frac{6}{35}h~~&\frac{6}{7}h-\frac{6}{35}\\
      \end{matrix}\right),$$
$$W(u)=\left(\begin{matrix}
            \frac{2}{5}u+\frac{8}{5}u^{3}\\
             \frac{12}{35}u+\frac{18}{35}u^{3}+\frac{8}{7}u^{5}\\
      \end{matrix}\right)\cdot\sigma(u).$$

\vskip 0.1 true cm
\noindent{\bf Proof.} By direct computations, we have
$$
J_{i,2j+1}^{'}=2(2j+1)\int_{u(h)^{2}}^{\tilde{x}(h)}x^{i}y^{2j}\frac{\partial{y}}{\partial{h}}dx+2\tilde{x}(h)^{i}y(\tilde{x},h)^{2j+1}\frac{\partial{\tilde{x}}}{\partial{h}}
-2u(h)^{2i+2j+1}\frac{\partial{u(h)^{2}}}{\partial{h}}.\eqno(3.16)$$
Differential $H(x,0)=h$ and $H(u^{2},u)=h$ with respect to $h$, we have
$$\frac{\partial{x}}{\partial{h}}=\frac{1}{x(1-x)},~~\frac{\partial{u}}{\partial{h}}=\frac{1}{u+2u^{3}-2u^{5}}.$$
Through the analysis of the singular points of the system, we have $\left|\frac{\partial{\tilde{x}}}{\partial{h}}(h)\right|<\infty$. Notice $\frac{\partial{y}}{\partial{h}}=1/y$, we can obtain that
$$J_{i,2j+1}^{'}=(2j+1)J_{i,2j-1}-4\sigma(u)u^{2i+2j+1}.$$
Similarly, we have
$$I_{i,2j+1}^{'}=(2j+1)I_{i,2j-1}+4\sigma(u)u^{2i+2j+1}.$$
Thus we have
$$J_{i,2j+1}=\frac{1}{2j+3}\left(J_{i,2j+3}^{'}+4\sigma(u)u^{2i+2j+3}\right),\eqno(3.17)$$
$$I_{i,2j+1}=\frac{1}{2j+3}\left(I_{i,2j+3}^{'}-4\sigma(u)u^{2i+2j+3}\right).\eqno(3.18)$$
Combining (3.8)--(3.13), we can obtain (3.14) and (3.15).
\vskip 0.3 true cm

\noindent{\bf Lemma 3.3.} Let $J_{0}=J_{0,1}+I_{0,1}$, $J_{1}=J_{1,1}+I_{1,1}$ and $D(h)=h(6h-1)$. Then $J_{0}$ and $J_{1}$ satisfy
$$J_{0}^{'}=\frac{1}{D(h)}\left[k_{0,0}(h)J_{0}+k_{0,1}(h)J_{1}\right],\eqno(3.19)$$
$$J_{1}^{'}=\frac{1}{D(h)}\left[k_{1,0}(h)J_{0}+k_{1,1}(h)J_{1}\right],\eqno(3.20)$$
where $k_{0,0}(h)=5h-1$, $k_{0,1}(h)=6/7$, $k_{1,0}(h)=-h$ and $k_{1,1}(h)=7h$.

\vskip 0.1 true cm
\noindent{\bf Proof.} By (3.14), we have
$${\rm det}(B_{1}h+C_{1}){\bf V^{'}_{1}(h)}=(B_{1}h+C_{1})^{*}\left({\bf V_{1}(h)}-W(u)\right). \eqno(3.21)$$
Combining $\frac{1}{2}u^{2}+\frac{1}{2}u^{4}-\frac{1}{3}u^{6}=h$, we can obtain
$$J_{0,1}^{'}=\frac{1}{D(h)}\left[k_{0,0}(h)J_{0,1}+k_{0,1}(h)J_{1,1}+w_{1}(u)\right],$$
$$J_{1,1}^{'}=\frac{1}{D(h)}\left[k_{1,0}(h)J_{0,1}+k_{1,1}(h)J_{1,1}+w_{2}(u)\right],$$
where
$$w_{1}(u)=-\sigma(u)\left(\frac{38}{35}u^{5}+\frac{4}{7}u^{7}-\frac{16}{35}u^{9}\right),$$
$$w_{2}(u)=-\sigma(u)\left(\frac{6}{35}u^{3}+\frac{12}{35}u^{5}+\frac{26}{35}u^{7}+\frac{4}{7}u^{9}-\frac{16}{35}u^{11}\right).$$
Similarly, we have
$$I_{0,1}^{'}=\frac{1}{D(h)}\left[k_{0,0}(h)I_{0,1}+k_{0,1}(h)I_{1,1}-w_{1}(u)\right],$$
$$I_{1,1}^{'}=\frac{1}{D(h)}\left[k_{1,0}(h)I_{0,1}+k_{1,1}(h)I_{1,1}-w_{2}(u)\right].$$
Therefore, we can get (3.19) and (3.20) by $J_{0}=J_{0,1}+I_{0,1}$ and $J_{1}=J_{1,1}+I_{1,1}$. This ends the proof.

\vskip 0.3 true cm

\noindent{\bf Lemma 3.4.} Let $\phi_{1}(h)=\alpha(h)J_{0}+\beta(h)J_{1}$. Then for $h\in(0,\frac{1}{6})$, there exist polynomials $P_{i}(h)(i=0,1,2)$ such that $L(h)\phi_{1}(h)=0$, where
$$L(h)=P_{2}(h)D(h)\frac{d^{2}}{dh^{2}}+P_{1}(h)D(h)\frac{d}{dh}+P_{0}(h),\eqno(3.22)$$
where
${\rm deg}P_{2}(h)\leq n_{1}$, ${\rm deg}P_{1}(h)\leq n_{1}-1$, ${\rm deg}P_{0}(h)\leq n_{1}-2$,
and
$$n_{1}=\left[(n-1)/2\right]+\left[n/2\right]+2.$$
In addition, we have
$$M_{1}(h):=L(h)M(h)={\tilde \gamma}(h)I_{0,1}+{\tilde \eta}(h)I_{1,1}+{\tilde \Phi}(u),\eqno(3.23)$$
where
${\rm deg}{\tilde \gamma}(h)\leq 2\left[\frac{n-1}{2}\right]+\left[\frac{n}{2}\right]+2$, ${\rm deg}{\tilde \eta}(h)\leq \left[\frac{n-1}{2}\right]+2\left[\frac{n}{2}\right]+1$ and ${\tilde \Phi}(u)$ is a rational fraction of $u$.

\vskip 0.1 true cm
\noindent{\bf Proof.} Differential both sides (3.21) and combining (3.14), we have
$$J_{0,1}^{''}=\frac{1}{D(h)}\left(-\frac{5}{6}J_{0,1}+w_{1}^{*}(u)\right),$$
$$J_{1,1}^{''}=\frac{1}{D(h)}\left(-J_{0,1}+\frac{7}{6}J_{1,1}+w_{2}^{*}(u)\right),$$
where
$$w_{1}^{*}(u)=-\frac{2}{35}u(7-39u^{2}-72u^{4}-46u^{6}-24u^{8}+32u^{10})\sigma(u)^{3},$$
$$w_{2}^{*}(u)=-\frac{2}{35}u(3+3u^{2}+77u^{4}+176u^{6}+22u^{8}-120u^{10}+32u^{12})\sigma(u)^{3}.$$
Similarly, we have
$$I_{0,1}^{''}=\frac{1}{D(h)}\left(\frac{5}{6}I_{0,1}-w_{1}^{*}(u)\right),$$
$$I_{1,1}^{''}=\frac{1}{D(h)}\left(-I_{0,1}+\frac{7}{6}I_{1,1}-w_{2}^{*}(u)\right).$$
Hence,
$$J_{0}^{''}=-\frac{5}{6D(h)}J_{0},~~J_{1}^{''}=\frac{1}{D(h)}\left(-J_{0}+\frac{7}{6}J_{1}\right).\eqno(3.24)$$
Suppose that
$$P_{2}(h)=\sum_{k=0}^{s}p_{2,k}h^{k},~~P_{1}(h)=\sum_{k=0}^{s-1}p_{1,k}h^{k},~~P_{0}(h)=\sum_{k=0}^{s-2}p_{0,k}h^{k}.$$
Using (3.19), (3.20) and (3.24), we can obtain
$$\begin{aligned}
L(h)\phi_{1}(h)&=P_{2}(h)D(h)\phi_{1}^{''}(h)+P_{1}(h)D(h)\phi_{1}^{'}(h)+P_{0}(h)\phi_{1}(h)\\
&=X(h)J_{0}+Y(h)J_{1},
\end{aligned}$$
where $X(h)$ and $Y(h)$ are polynomials with degree no more than $2\left[(n-1)/2\right]+\left[n/2\right]+2$ and $\left[(n-1)/2\right]+2\left[n/2\right]+1$ respectively. Let
$$X(h)=\sum_{i=0}^{{\rm deg}X}x_{i}h^{i},~~Y(h)=\sum_{j=0}^{{\rm deg}Y}y_{j}h^{j},$$
$x_{i}$ and $y_{j}$ are expressed by $p_{2,k}$, $p_{1,k}$ linearly. So $L(h)\phi_{1}(h)=0$ is satisfied if we let
$$x_{i}=0,~~y_{j}=0,~~(0\leq i\leq {\rm deg}X,~0\leq j\leq {\rm deg}Y).\eqno(3.25)$$
System (3.25) is a homogeneous linear equation with $3\left[(n-1)/2\right]+3\left[n/2\right]+5$ equations about $3\left[(n-1)/2\right]+3\left[n/2\right]+6$ variables of $p_{2,k}$, $p_{1,k}$ and $p_{0,k}$. It follows that from the theory of linear algebra that exist $p_{2,k}$, $p_{1,k}$ and $p_{0,k}$ such that the result holds. Let
$$\Phi_{1}(u)=\Phi^{'}(u)D(u)\frac{\partial{u}}{\partial{h}}+w_{1}(u)\left(\alpha(u)-\gamma(u)\right)+w_{2}(u)\left(\beta(u)-\eta(u)\right),\eqno(3.26)$$
$$\Phi_{2}(u)=\left(\Phi^{'}(u)\frac{\partial{u}}{\partial{h}}\right)^{'}_{u}\frac{\partial{u}}{\partial{h}}D(u)+2w_{1}(\alpha^{'}-\gamma^{'})+2w_{2}(\beta^{'}-\eta^{'})
+w_{1}^{*}(\alpha-\gamma)+w_{2}^{*}(\beta-\eta).\eqno(3.27)$$
By the same progress, we can obtain
$$\begin{aligned}
L(h)M(h)&=P_{2}(h)D(h)M^{''}(h)+P_{1}(h)D(h)M^{'}(h)+P_{0}(h)M(h)\\
&=X(h)J_{0,1}+Y(h)J_{1,1}+{\tilde \gamma}(h)I_{0,1}+{\tilde \eta}(h)I_{1,1}+{\tilde \Phi}(u),
\end{aligned}$$
where
$${\rm deg}{\tilde \gamma}(h)\leq 2\left[\frac{n-1}{2}\right]+\left[\frac{n}{2}\right]+2,~~{\rm deg}{\tilde \eta}(h)\leq \left[\frac{n-1}{2}\right]+2\left[\frac{n}{2}\right]+1,$$
and $${\tilde \Phi}(u)=P_{2}(u)\Phi_{2}(u)+P_{1}(u)\Phi_{1}(u)+P_{0}(u)\Phi(u).\eqno(3.28)$$
This ends the proof.
\vskip 0.3 true cm
Let
$${\tilde \Phi}_{1}(u)={\tilde \Phi}^{'}(u)D(u)\frac{\partial{u}}{\partial{h}}-w_{1}(u){\tilde \gamma}(u)-w_{2}(u){\tilde \eta}(u),\eqno(3.29)$$

$${\tilde \Phi}_{2}(u)=\left({\tilde \Phi}^{'}(u)\frac{\partial{u}}{\partial{h}}\right)^{'}_{u}\frac{\partial{u}}{\partial{h}}D(u)-2w_{1}{\tilde \gamma}^{'}-2w_{2}{\tilde \eta}^{'}-w_{1}^{*}{\tilde \gamma}-w_{2}^{*}{\tilde \eta}.\eqno(3.30)$$

Similar to the proof of Lemma 3.4, we can obtain the following Lemma.
\vskip 0.3 true cm

\noindent{\bf Lemma 3.5.} Let $\phi_{2}(h)={\tilde \gamma}(h)J_{0}+{\tilde \eta}(h)J_{1}$. Then for $h\in(0,\frac{1}{6})$, there exist polynomials ${\tilde P}_{i}(h)(i=0,1,2)$ such that ${\tilde L}(h)\phi_{2}(h)=0$, where
$${\tilde L}(h)={\tilde P}_{2}(h)D(h)\frac{d^{2}}{dh^{2}}+{\tilde P}_{1}(h)D(h)\frac{d}{dh}+{\tilde P}_{0}(h),\eqno(3.31)$$
and
$$M_{2}(h):={\tilde L}(h)M_{1}(h)={\hat \Phi}(u),\eqno(3.32)$$
where
${\rm deg}{\tilde P}_{2}(h)\leq n_{2}$, ${\rm deg}{\tilde P}_{1}(h)\leq n_{2}-1$, ${\rm deg}{\tilde P}_{0}(h)\leq n_{2}-2$,
$${\hat \Phi}(u)={\tilde P}_{2}(u){\tilde \Phi}_{2}(u)+{\tilde P}_{1}(u){\tilde \Phi}_{1}(u)+{\tilde P}_{0}(u){\tilde \Phi}(u),\eqno(3.33)$$
and
$$n_{2}=3\left(\left[(n-1)/2\right]+\left[n/2\right]\right)+6$$.

\vskip 0.3 true cm

\noindent{\bf Proof for the Theorem 1.1.} Denote $\mathcal{L}(u_{1},~u_{2},~u_{3},...,u_{n-1},~u_{n})$ is the linear combination of function $u_{1},~u_{2},~u_{3},...,u_{n-1},~u_{n}$. First, we will analysis the structure of ${\tilde \Phi}(u)$ and ${\hat \Phi}(u)$. Through direct computation we can get for polynomial $f(h)$ with degree no more than $n_{0}$, it can become the linear combination of $1,~u^{2},~u^{4},...,u^{2n-2},~u^{2n}~(n\leq 3n_{0})$ by $h=\frac{1}{2}u^{2}+\frac{1}{2}u^{4}-\frac{1}{3}u^{6}$. Thus, from the expression of $\Phi(u),~w_{i}(u),~w_{i}^{*}(u)~(i=1,2)$, we can obtain
$$\Phi_{1}(u)=\sigma(u)\cdot\mathcal{L}(u,~u^{3},~u^{5},...,u^{6\left[\frac{n-1}{2}\right]+11},~u^{6\left[\frac{n-1}{2}\right]+13}),$$
$$\Phi_{2}(u)=\sigma(u)^{3}\cdot\mathcal{L}(1,~u^{2},~u^{4},...,u^{6\left[\frac{n-1}{2}\right]+14},~u^{6\left[\frac{n-1}{2}\right]+16})/u.$$
Therefore, we have
$$\begin{aligned}
{\tilde \Phi}(u)&=P_{2}(h)\Phi_{2}(u)+P_{1}(h)\Phi_{1}(u)+P_{0}(h)\Phi(u)\\
&=\sigma(u)^{3}\cdot\mathcal{L}(1,~u^{2},~u^{4},...,u^{2k_{1}-2},~u^{2k_{1}})/u,
\end{aligned}$$
with $k_{1}=6\left[\frac{n-1}{2}\right]+3\left[\frac{n}{2}\right]+14$.
Similarly, we can get
$$\begin{aligned}
{\hat \Phi}(u)&={\tilde P}_{2}(u(h)){\tilde \Phi}_{2}(u)+{\tilde P}_{1}(u(h)){\tilde \Phi}_{1}(u)+{\tilde P}_{0}(u(h)){\tilde \Phi}(u)\\
&=\sigma(u)^{7}\cdot\mathcal{L}(1,~u^{2},~u^{4},...,u^{2k_{2}-2},~u^{2k_{2}})/u^{3},
\end{aligned}$$
with $k_{2}=15\left[\frac{n-1}{2}\right]+12\left[\frac{n}{2}\right]+41$.
Suppose $u_{*}$ satisfys $\frac{1}{2}u_{*}^{2}+\frac{1}{2}u_{*}^{4}-\frac{1}{3}u_{*}^{6}=\frac{1}{6}$, then the number of zeros of function $f(h)$ are same with the number of zeros of function $f(u)$ in $h\in(0,\frac{1}{6})$ and $u\in(0,u_{*})$. Notice $1+2u^2-2u^4\neq 0$ and function $1,~u,~u^{2},...,u^{i}$ are ECT-system in $u\in(0,u_{*})$, we have
$$\#\left\{{\tilde \Phi}(u)=0,u\in(0,u_{*})\right\}\leq k_{1},~~\#\left\{{\hat \Phi}(u)=0,u\in(0,u_{*})\right\}\leq k_{2}.$$

Next, we analysis the number of zeros of function $\phi_{1}(h)$ and $\phi_{2}(h)$. Since
$$J_{0}(h)=J_{0,1}(h)+I_{0,1}(h)=\oint_{L_{h}^{+}\cup L_{h}^{-}}ydx=\iint\limits_{L_{h}^{+}\cup L_{h}^{-}}dxdy\neq 0,$$
we let
$$S(h)=J_{1}(h)/J_{0}(h),~~F(h)=\frac{\phi_{1}(h)}{J_{0}(h)}=\alpha(h)+\beta(h)S(h).$$
Then $F(h)$ satisfy the following Riccati equation
$$D(h)\beta(h)F^{'}(h)=-k_{0,1}\beta(h)F(h)^{2}+N_{1}(h)F(h)+N_{2}(h),$$
with $N_{1}(h)$ and $N_{2}(h)$ are polynomials with degree no more than $\left[\frac{n-1}{2}\right]+\left[\frac{n}{2}\right]-1$ and $2\left[\frac{n-1}{2}\right]+\left[\frac{n}{2}\right]-1$ respectively. Thus we have
$$\begin{aligned}
\#\left\{\phi_{1}(h)=0,h\in(0,1/6)\right\}&\leq \#\left\{D(h)\beta(h)=0,h\in(0,1/6)\right\}+\#\left\{N_{2}(h)=0,h\in(0,1/6)\right\}+1\\
&\leq 2\left(\left[(n-1)/2\right]+\left[n/2\right]\right)+1.
\end{aligned}$$

Similarly, we get
$$\#\left\{\phi_{2}(h)=0,h\in(0,1/6)\right\}\leq 6\left(\left[(n-1)/2\right]+\left[n/2\right]\right)+9.$$

Therefore By Definition 2.2 and Lemma 2.3--2.4, we can obtain
$$\begin{aligned}
\#\left\{M_{1}(h)=0,\right.&\left.h\in(0,1/6)\right\}\leq \#\left\{{\hat \Phi}(u)=0,u\in(0,u_{*})\right\}\\
&+3\left(\#\left\{{\tilde P}_{2}(h)D(h)=0,h\in(0,1/6)\right\}+\#\left\{\phi_{2}(h)=0,h\in(0,1/6)\right\}\right)+2\\
&\leq 42\left[(n-1)/2\right]+39\left[n/2\right]+94,
\end{aligned}
$$
$$\begin{aligned}
\#\left\{M(h)=0,\right.&\left.h\in(0,1/6)\right\}\leq \#\left\{M_{1}(h)=0,h\in(0,1/6)\right\}\\
&+3\left(\#\left\{P_{2}(h)D(h)=0,h\in(0,1/6)\right\}
+\#\left\{\phi_{1}(h)=0,h\in(0,1/6)\right\}\right)+2\\
&\leq 51\left[(n-1)/2\right]+48\left[n/2\right]+111,
\end{aligned}
$$

\vskip 0.5 true cm
\centerline{\bf{ $\S$4}. Proof of Theorem 1.2}
\vskip 0.3 true cm

This section we consider system (1.3). Similar to the Lemma 2.1, we can obtain
$$\begin{aligned}
M(h)&=\int_{L_{h}^{+}}q^{+}(x,y)dx-p^{+}(x,y)dy+\int_{L_{h}^{-}}q^{-}(x,y)-p^{-}(x,y)dy\\
&+\int_{\tilde{L}_{h}^{+}}\tilde{q}^{+}(x,y)dx-\tilde{p}^{+}(x,y)dy+\int_{\tilde{L}_{h}^{-}}\tilde{q}^{-}(x,y)-\tilde{p}^{-}(x,y)dy,
\end{aligned}\eqno(4.1)$$
where
$$L_{h}^{+}(L_{h}^{-})=\{(x,y)|H(x,y)=h,x>0,y>0(x<0,y>0)\},$$
$$\tilde{L}_{h}^{+}(\tilde{L}_{h}^{-})=\{(x,y)|H(x,y)=h,x>0,y<0(x<0,y<0)\}.$$
Let
$$J_{i,j}(h)=\int_{{L}_{h}^{+}}x^{i}y^{j}dx,~~I_{i,j}(h)=\int_{L_{h}^{-}}x^{i}y^{j}dx,$$
$$\tilde{J}_{i,j}(h)=\int_{{\tilde{L}}_{h}^{+}}x^{i}y^{j}dx,~~\tilde{I}_{i,j}(h)=\int_{\tilde{L}_{h}^{-}}x^{i}y^{j}dx.$$
Then we have $\tilde{J}_{i,j}(h)=(-1)^{j+1}J_{i,j}(h)$, $\tilde{I}_{i,j}(h)=(-1)^{j+1}I_{i,j}(h)$. Then we have
$$\begin{aligned}
M(h)&=(\rho_{0,0}^{+}-\tilde{\rho}_{0,0}^{+})J_{0,0}+(\rho_{1,0}^{+}-\tilde{\rho}_{1,0}^{+})J_{1,0}+(\rho_{0,1}^{+}+\tilde{\rho}_{0,1}^{+})J_{0,1}\\
&+(\rho_{2,0}^{+}-\tilde{\rho}_{2,0}^{+})J_{2,0}+(\rho_{1,1}^{+}+\tilde{\rho}_{1,1}^{+})J_{1,1}+(\rho_{0,2}^{+}-\tilde{\rho}_{0,2}^{+})J_{0,2}\\
&+(\rho_{0,0}^{-}-\tilde{\rho}_{0,0}^{-})I_{0,0}+(\rho_{1,0}^{-}-\tilde{\rho}_{1,0}^{-})I_{1,0}+(\rho_{0,1}^{-}+\tilde{\rho}_{0,1}^{-})I_{0,1}\\
&+(\rho_{2,0}^{-}-\tilde{\rho}_{2,0}^{-})I_{2,0}+(\rho_{1,1}^{-}+\tilde{\rho}_{1,1}^{-})I_{1,1}+(\rho_{0,2}^{-}-\tilde{\rho}_{0,2}^{-})I_{0,2}\\
&:=p_{0,0}J_{0,0}+p_{1,0}J_{1,0}+p_{0,1}J_{0,1}+p_{2,0}J_{2,0}+p_{1,1}J_{1,1}+p_{0,1}J_{0,1}\\
&+q_{0,0}I_{0,0}+q_{1,0}I_{1,0}+q_{0,1}I_{0,1}+q_{2,0}I_{2,0}+q_{1,1}I_{1,1}+q_{0,1}I_{0,1}
\end{aligned}\eqno(4.2)$$
where $\rho_{i,j}^{\pm}(\tilde{\rho}_{i,j}^{\pm})=q_{i,j}^{\pm}(\tilde{q}_{i,j}^{\pm})+\frac{i+1}{j}p_{i+1,j-1}^{\pm}(\tilde{p}_{i+1,j-1}^{\pm})~(j\geq 1)$ and $\rho_{i,0}^{\pm}(\tilde{\rho}_{i,0}^{\pm})=p_{i,0}^{\pm}(\tilde{p}_{i,0}^{\pm})$. By direct computation, we have
$$J_{0,0}=\tilde{x}(h),~~I_{0,0}=-\bar{x}(h),~~J_{1,0}=\frac{1}{2}\tilde{x}(h)^{2},~~I_{1,0}=-\frac{1}{2}\bar{x}(h)^{2},~~J_{2,0}=\frac{1}{3}\tilde{x}(h)^{3}.\eqno(4.3)$$
$$I_{2,0}=-\frac{1}{3}\bar{x}(h)^{3},~~J_{0,2}=2h\tilde{x}(h)-\frac{1}{3}\tilde{x}(h)^{3}+\frac{1}{6}\tilde{x}(h)^{4},~~I_{0,2}=-2h\bar{x}(h)+\frac{1}{3}\bar{x}(h)^{3}-\frac{1}{6}\bar{x}(h)^{4},\eqno(4.4)$$
where
$$\tilde{x}(h)=\cos\left(\frac{\pi}{3}+\frac{1}{3}\arccos(12h-1)\right)+\frac{1}{2},~~\bar{x}(h)=-\cos\left(\frac{1}{3}\arccos(12h-1)\right)+\frac{1}{2}.$$

Let $u=x(\frac{1}{2}-\frac{1}{3}x)^{1/2}$, then for $x>0$ small enough, by the implicit theorem we have
$$\begin{aligned}
x=\phi(u)&=-\sqrt{2}u-\frac{\sqrt{2}}{3}u^2-\frac{11}{54}\sqrt{2}u^3-\frac{\sqrt{2}}{8}u^4-\frac{379}{3240}\sqrt{2}u^5-\frac{565}{5832}\sqrt{2}u^6-\frac{751}{9072}\sqrt{2}u^7\\
&-\frac{1687}{23328}\sqrt{2}u^8-\frac{161809}{2519424}\sqrt{2}u^9-\frac{727783}{12597120}\sqrt{2}u^{10}-\frac{8730965}{166281984}\sqrt{2}u^{11}+o(u^{11}).
\end{aligned}\eqno(4.5)$$
Therefore
$$\begin{aligned}
J_{0,0}&=\int_{0}^{\tilde{x}(h)}\sqrt{2h+\frac{2}{3}x^3-x^2}dx=\int_{0}^{h^{\frac{1}{2}}}\sqrt{2h-2u^2}\phi^{'}(u)du\\
&=\int_{0}^{h^{\frac{1}{2}}}\sqrt{2h-2(h^\frac{1}{2}u)^{2}}\phi^{'}(h^{\frac{1}{2}}u)d(h^{\frac{1}{2}}u)=\sqrt{2}h\int_{0}^{1}\sqrt{1-2u^{2}}\phi^{'}(h^{\frac{1}{2}}u)d(u)\\
&=\sqrt{2}h\left(-\frac{1}{4}\sqrt{2}\pi-\frac{2}{9}\sqrt{2}h^{\frac{1}{2}}-\frac{11}{288}\sqrt{2}\pi h-\frac{1}{15}\sqrt{2}h^{\frac{3}{2}}-\frac{379}{20736}\sqrt{2}\pi h^2-\frac{226}{5103}\sqrt{2}h^\frac{5}{2}\right.\\
&\left.\quad\quad\quad-\frac{3755}{331776}\sqrt{2}\pi h^3
-\frac{964}{32805}\sqrt{2}h^{\frac{7}{2}}-\frac{1132663}{143327232}\sqrt{2}\pi h^4-\frac{207938}{9743085}\sqrt{2}h^{\frac{9}{2}}\right.\\
&\left.\quad\quad\quad-\frac{61116755}{10319560704}\sqrt{2}\pi h^5+o(h^5)\right),
\end{aligned}\eqno(4.6)$$
$$\begin{aligned}
I_{0,0}&=\int_{\bar{x}(h)}^{0}\sqrt{2h+\frac{2}{3}x^3-x^2}dx=\int_{-h^{\frac{1}{2}}}^{0}\sqrt{2h-2u^2}\phi^{'}(u)du=\sqrt{2}h\int_{-1}^{0}\sqrt{1-u^{2}}\phi^{'}(h^{\frac{1}{2}}u)du\\
&=\sqrt{2}h\left(-\frac{1}{4}\sqrt{2}\pi+\frac{2}{9}\sqrt{2}h^{\frac{1}{2}}-\frac{11}{288}\sqrt{2}\pi h+\frac{1}{15}\sqrt{2}h^{\frac{3}{2}}-\frac{379}{20736}\sqrt{2}\pi h^2+\frac{226}{5103}\sqrt{2}h^{\frac{5}{2}}\right.\\
&\left.\quad\quad\quad-\frac{3755}{331776}\sqrt{2}\pi h^3+\frac{964}{32805}\sqrt{2}h^{\frac{7}{2}}-\frac{1132663}{143327232}\sqrt{2}\pi h^4+\frac{207938}{9743085}\sqrt{2}h^{\frac{9}{2}}\right.\\
&\left.\quad\quad\quad-\frac{61116755}{10319560704}\sqrt{2}\pi h^5+o(h^5)\right).
\end{aligned}\eqno(4.7)$$
Similarly, we can obtain
$$\begin{aligned}
J_{1,1}(I_{1,1})&=\sqrt{2}h\left(\pm\frac{2}{3}h^{\frac{1}{2}}+\frac{1}{8}\pi h\pm\frac{112}{405}h^{\frac{3}{2}}+\frac{625}{10368}\pi h^2\pm\frac{20924}{127575}h^{\frac{5}{2}}\right.\\
&\left.+\frac{65863}{1492992}\pi h^3
\pm\frac{319442}{2679075}h^{\frac{7}{2}}+\frac{5919829}{179159040}\pi h^4\pm\frac{218941144}{2387055825}h^{\frac{9}{2}}\right.\\
&\left.+\frac{336369143}{12899450880}\pi h^5+o(h^{5})\right).
\end{aligned}\eqno(4.8)$$

In addition, for $h>0$ sufficiently small
$$\begin{aligned}
\tilde{x}(h)&=\sqrt{2}h^{\frac{1}{2}}+\frac{2}{3}h+\frac{5}{9}\sqrt{2}h^{\frac{3}{2}}+\frac{32}{27}h^2+\frac{77}{54}\sqrt{2}h^{\frac{5}{2}}+\frac{896}{243}h^3+\frac{2431}{486}\sqrt{2}h^{\frac{7}{2}}+\frac{10240}{729}h^4\\
&+\frac{1062347}{52488}\sqrt{2}h^{\frac{9}{2}}+\frac{1171456}{19683}h^5+\frac{14003665}{157464}\sqrt{2}h^\frac{11}{2}+\frac{47710208}{177147}h^6+o(h^6),
\end{aligned}\eqno(4.9)$$
$$\begin{aligned}
\bar{x}(h)&=-\sqrt{2}h^{\frac{1}{2}}+\frac{2}{3}h-\frac{5}{9}\sqrt{2}h^{\frac{3}{2}}+\frac{32}{27}h^2-\frac{77}{54}\sqrt{2}h^{\frac{5}{2}}+\frac{896}{243}h^3-\frac{2431}{486}\sqrt{2}h^{\frac{7}{2}}+\frac{10240}{729}h^4\\
&-\frac{1062347}{52488}\sqrt{2}h^{\frac{9}{2}}+\frac{1171456}{19683}h^5-\frac{14003665}{157464}\sqrt{2}h^\frac{11}{2}+\frac{47710208}{177147}h^6+o(h^6).
\end{aligned}\quad\eqno(4.10)$$
Therefore, combining (4.3)--(4.4) and (4.6)--(4.10), we can get
$$
M(h)=\sum\limits_{i\geq 1}\delta_{i}h^{\frac{i}{2}},\eqno(4.11)
$$
where
$$\begin{aligned}
&\delta_{1}=\sqrt{2}(p_{0,0}+q_{0,0}),\\
&\delta_{2}=-\frac{2}{3}(q_{0,0}-p_{0,0})-\frac{\pi}{2}(q_{0,1}+p_{0,1})-(q_{1,0}-p_{1,0}),\\
&\delta_{3}=\frac{5}{9}\sqrt{2}(q_{0,0}+p_{0,0})+\frac{4}{9}(q_{0,1}-p_{0,1})-\frac{4}{3}\sqrt{2}(q_{0,2}-p_{0,2})+\frac{2}{3}\sqrt{2}(q_{1,0}+p_{1,0})\\
&\qquad-\frac{2}{3}\sqrt{2}(q_{1,1}-p_{1,1})+\frac{2}{3}\sqrt{2}(q_{2,0}+p_{2,0}),\\
&\delta_{4}=-\frac{32}{27}(q_{0,0}-p_{0,0})-\frac{11}{144}\pi(q_{0,1}+p_{0,1})+\frac{2}{3}(q_{0,2}+p_{0,2})-\frac{4}{3}(q_{1,0}+p_{1,0})\\
&\qquad+\frac{1}{8}\sqrt{2}\pi(q_{1,1}+p_{1,1})-\frac{4}{3}(q_{2,0}-p_{2,0}),\\
&\delta_{5}=\frac{77}{54}\sqrt{2}(q_{0,0}+p_{0,0})+\frac{2}{15}(q_{0,1}-p_{0,1})+\frac{4}{9}\sqrt{2}(q_{0,2}+p_{0,2})+\frac{14}{9}\sqrt{2}(q_{1,0}+p_{1,0})\\
&\qquad-\frac{112}{405}\sqrt{2}(q_{1,1}-p_{1,1})+\frac{14}{9}\sqrt{2}(q_{2,0}+p_{2,0}),\\
&\delta_{6}=-\frac{896}{243}(q_{0,0}-p_{0,0})-\frac{379}{10368}\pi(q_{0,1}+p_{0,1})-\frac{64}{81}(q_{0,2}-p_{0,2})-\frac{320}{81}(q_{1,0}-p_{1,0})\\
&\qquad+\frac{625}{10368}\sqrt{2}\pi(q_{1,1}+p_{1,1})-\frac{320}{81}(q_{2,0}-p_{2,0}),\\
\end{aligned}$$
$$\begin{aligned}
&\delta_{7}=\frac{2431}{486}\sqrt{2}(q_{0,0}+p_{0,0})+\frac{452}{5103}(q_{0,1}-p_{0,1})+\frac{22}{27}\sqrt{2}(q_{0,2}+p_{0,2})+\frac{143}{27}\sqrt{2}(q_{1,0}+p_{1,0})\\
&\qquad-\frac{20924}{127575}\sqrt{2}(q_{1,1}-p_{1,1})+\frac{143}{27}\sqrt{2}(q_{2,0}+p_{2,0}),\\
&\delta_{8}=-\frac{10240}{729}(q_{0,0}-p_{0,0})-\frac{3755}{165888}\pi(q_{0,1}+p_{0,1})-\frac{448}{243}(q_{0,2}-p_{0,2})-\frac{3584}{243}(q_{1,0}-p_{1,0})\\
&\qquad+\frac{65863}{1492992}\sqrt{2}\pi(q_{1,1}+p_{1,1})-\frac{3584}{243}(q_{2,0}-p_{2,0}),\\
&\delta_{9}=\frac{1062347}{52488}\sqrt{2}(q_{0,0}+p_{0,0})+\frac{1928}{32805}(q_{0,1}-p_{0,1})+\frac{4862}{2187}\sqrt{2}(q_{0,2}+p_{0,2})\\
&\qquad+\frac{46189}{2187}\sqrt{2}(q_{1,0}+p_{1,0})-\frac{319442}{2679075}\sqrt{2}(q_{1,1}-p_{1,1})+\frac{46189}{2187}\sqrt{2}(q_{2,0}+p_{2,0}),\\
&\delta_{10}=-\frac{1171456}{19683}(q_{0,0}-p_{0,0})-\frac{1132663}{71663616}\pi(q_{0,1}+p_{0,1})-\frac{4096}{729}(q_{0,2}-p_{0,2})\\
&\qquad-\frac{45056}{729}(q_{1,0}-p_{1,0})+\frac{5919829}{179159040}\sqrt{2}\pi(q_{1,1}+p_{1,1})-\frac{45056}{729}(q_{2,0}-p_{2,0}),\\
&\delta_{11}=\frac{14003665}{157464}\sqrt{2}(q_{0,0}+p_{0,0})+\frac{415876}{9743085}(q_{0,1}-p_{0,1})+\frac{96577}{13122}\sqrt{2}(q_{0,2}+p_{0,2})\\
&\qquad+\frac{2414425}{26244}\sqrt{2}(q_{1,0}+p_{1,0})-\frac{218941144}{2387055825}\sqrt{2}(q_{1,1}-p_{1,1})+\frac{2414425}{26244}\sqrt{2}(q_{2,0}+p_{2,0}),\\
&\delta_{12}=-\frac{47710208}{177147}(q_{0,0}-p_{0,0})-\frac{61116755}{5159780352}\pi(q_{0,1}+p_{0,1})-\frac{1171456}{59049}(q_{0,2}-p_{0,2})\\
&\qquad-\frac{16400384}{59049}(q_{1,0}-p_{1,0})+\frac{336369143}{12899450880}\sqrt{2}\pi(q_{1,1}+p_{1,1})-\frac{16400384}{59049}(q_{2,0}-p_{2,0}),\\
\end{aligned}$$
and
$$\delta_{i}=\mathcal{L}(\delta_{1},\delta_{2},\delta_{3},...,\delta_{12}),~i\geq 13.\eqno(4.12)$$
Then we calculate that
$$
\qquad\frac{\partial\left(\delta_{1},\delta_{2},\delta_{3},\delta_{4},\delta_{5},\delta_{6},\delta_{7},\delta_{8},\delta_{9},\delta_{10},\delta_{11},\delta_{12}\right)}
{\partial\left(p_{0,0},p_{1,0},p_{0,1},p_{2,0},p_{1,1},p_{0,2},q_{0,0},q_{1,0},q_{0,1},q_{2,0},q_{1,1},q_{0,2}\right)}$$
$$=\frac{59886739950651665292703225049}{15759296625811548028684506048}\pi^{2},\eqno(4.13)$$
which means that $\delta_{1},\delta_{2},...,\delta_{12}$ can be chosen arbitrary. On the other hands, by (4.12) we have
$$M(h)=h^{\frac{1}{2}}\delta_{1}\left(1+P_{1}(h,\delta)\right)+h\delta_{2}\left(1+P_{2}(h,\delta)\right)+h^{\frac{3}{2}}\delta_{3}\left(1+P_{3}(h,\delta)\right)+\cdots+h^{6}\delta_{12}\left(1+P_{12}(h,\delta)\right),$$
where $\delta=\left(\delta_{1},\delta_{2},\delta_{3}...,\delta_{12}\right)$ and
$$P_{i}\in C^{\infty},~ P_{i}(0,\delta)=0,~(i=1,2,...,12).$$
Therefore, $M(h)\equiv0$ when $\delta_{1}=\delta_{2}=\cdots=\delta_{12}=0$. By [11], we can get the Theorem 1.2 combining (4.13). This ends the proof.

\vskip 0.5 true cm


\begin{thebibliography}{30}
\bibitem{ref1} {E. Barbashin, Introduction to the theory of stability, Noordhoff, Groningen, 2017.}
\bibitem{ref2} {J. Bastos, C. A. Buzzi, J. Llibre, D. D. Novaes, Melnikov analysis in nonsmooth differential systems with nonlinear switching manifold, J. Differential Equations  267 (2019) 3748-3767.}
\bibitem{ref3} {W. Cui, B. Li, Y. Zhang, Lowest upper bounds of number of limit cycles on Bogdanov-Takens system under piecewise n-degree polynomial perturbation, Journal of Tianjin Normal University (Natural Science Edition) 40 (2020) 1-7.}
\bibitem{ref4} {X. Cen, C. Liu, L. Yang, M. Zhang, Limit cycles by perturbing quadratic isochronous centers inside piecewise polynomial differential systems, J. Differential Equations 265 (2018) 6083-6126.}
\bibitem{ref5} {J. -P. Francoise, H. Ji, D. Xiao, J. Yu, Global dynamics of a piecewise smooth system for brain lactate metabolism, Qual. Theory Dyn. Syst. 18 (2019) 315-332.}
\bibitem{ref6} {A. Gasull, J. Torregrosa, X. Zhang, Piecewise linear differential systems with an algebraic line of separation, Electron. J. Differential Equations 19 (2020) 1-14.}
\bibitem{ref7} {L. Gavrilov, I. D. Iliev, Quadratic perturbations of quadratic codimension-four centers, J. Math. Anal. Appl. 357 (2009) 69-76.}
\bibitem{ref8} {M. Grau, F. Ma$\tilde{n}$osas, J. Villadelprat, A Chebyshev criterion for Abelian integrals, Trans. Am. Math. Soc. 363 (2011) 109-129.}
\bibitem{ref9} {E. Horozov, I. D. Iliev, On the number of limit cycles in perturbations of quadratic Hamiltonian systems, Proc. Lond. Math. Soc. 69 (1994) 198-224.}
\bibitem{ref10} {M. Han, Bifurcation Theory of Limit Cycles, Science Press, Beijing, 2013.}
\bibitem{ref11} {M. Han, L. Sheng, Bifurcation of limit cycles in piecewise smooth systems via Melnikov function, J. Appl. Anal. Comput. 5 (2015) 809-815.}
\bibitem{ref12} {T. Ito, A Filippov solution of a system of differential equations with discontinuous right-hand sides, Econom. Lett. 4 (1979) 349-354.}
\bibitem{ref13} {I. D. Iliev, High-order Melnikov functions for degenerate cubic Hamiltonians, Adv. Diff. Eqns. 1 (1996) 689-708.}
\bibitem{ref14} {V. Krivan, On the Gause predator-prey model with a refuge: a fresh look at the history, J. Theoret. Biol. 274 (2011) 67-73.}
\bibitem{ref15} {B. Li, Z. Zhang, A note of a G.S. Petrov's result about the weakened 16th Hilbert problem, J. Math. Anal. Appl. 190 (1995) 489-516.}
\bibitem{ref16} {S. Li, C. Liu, A linear estimate of the number of limit cycles for some planar piecewise smooth quadratic differential system, J. Math. Anal. Appl. 428 (2015) 1354-1367.}
\bibitem{ref17} {X. Liu, M. Han, Bifurcation of limit cycles by perturbing piecewise Hamiltonian systems, Int. J. Bifurcation and Chaos 20 (2010) 1379-1390.}
\bibitem{ref18} {J. Llibre, Y. Tang, Limit cycles of discontinuous piecewise quatratic and cubic polynomial perturbations of a linear center, Discrete Contin. Dyn. Syst. 24 (2019) 1769-1784.}
\bibitem{ref19} {J. Llibre, X. Zhang, Limit cycles for discontinuous planar piecewise linear differential systems separated by an algebraic curve, Int. J. Bifurcation and Chaos 29 (2019) 1950017.}
\bibitem{ref20} {M. Teixeira, Perturbation theory for non-smooth systems, in: Encyclopedia of Complexity and Systems Science, vol. 22, Springer, New York, 2009.}
\bibitem{ref21} {D. D. Novaes, J. Torregrosa, On extended Chebyshev systems with positive accuracy, J. Math. Anal. Appl. 448 (2017) 171-186.}
\bibitem{ref22} {O. Ramirez, A. M. Alves, Bifurcation of limit cycles by perturbing piecewise non-Hamiltonian systems with nonlinear switching manifold, Nonlinear Anal. Real World Appl. 57 (2021) 103188.}
\bibitem{ref23} {S. Sui, J. Yang, L. Zhao, On the number of limit cycles for generic Lotka-Volterra system and Bogdanov-Takens system under perturbations of piecewise smooth polynomials, Nonlinear Anal. Real World Appl. 49 (2019) 137-158.}
\bibitem{ref24} {H. Tian, M. Han, Limit cycle bifurcations of piecewise smooth near-Hamiltonian systems with switching curve, Discrete and Continuous Dynamical Systems Series B, 26 (2021) 5581-5599.}
\bibitem{ref25} {Y. Xiong, M. Han, Limit cycle bifurcations by perturbing a class of planar quintic vector fields, J. Differential Equations 269 (2020) 10964-10994.}
\bibitem{ref26} {J. Yang, L. Zhao, Bounding the number of limit cycles of discontinuous differential systems by using Picard-Fuchs equations, J. Differential Equations 264 (2018) 5734-5757.}
\bibitem{ref27} {J. Yang, Limit cycles appearing from the perturbation of differential systems with multiple switching curves, Chaos Solitons Fractals 135 (2020) 109764.}
\bibitem{ref28} {Q. Zhao, J. Yu, Poincar\`{e} maps of ``$<$'' -shape planar piecewise linear dynamical systems with a saddle, Int. J. Bifurcat. Chaos 29 (2019) 1590165.}
\bibitem{ref29} {C. Zou, C. Liu, J. Yang, On piecewise linear differential systems with $n$ limit cycles of arbitrary multiplicities in two zones, Qual. Theory Dyn. Syst. 18 (2019) 139-151.}
\end{thebibliography}
\end{document}